\documentclass{article}
\usepackage[utf8]{inputenc}
\usepackage[scale=0.75]{geometry}
\usepackage{amsmath}
\usepackage{mathrsfs}
\usepackage{txfonts}
\usepackage{amssymb} % equation align

\usepackage[title]{appendix}
\usepackage{graphicx}
\usepackage{booktabs}
\usepackage{tabularx}
\usepackage{booktabs}
\usepackage{multirow}
\usepackage{pgfplots}
\usepackage{subfigure}

\usepackage[ruled,linesnumbered]{algorithm2e}
\usepackage{setspace}

\usepackage[round, authoryear, sort&compress]{natbib}
\usepackage{hyperref}
\hypersetup{
    colorlinks=true,
    linkcolor=black,
    filecolor=gray,
    urlcolor=blue,
    citecolor=blue,
}

\RequirePackage[normalem]{ulem} %DIF PREAMBLE
\RequirePackage{color}\definecolor{RED}{rgb}{1,0,0}\definecolor{BLUE}{rgb}{0,0,1} %DIF PREAMBLE
\usepackage{authblk}
\usepackage{xpatch}
\xpatchcmd{\author}{\relax#1\relax}{\relax\detokenize{#1}\relax}{}{}

\title{A multi-period multi-product stochastic inventory problem with order-based loan}

\author[\empty]{Zhen Chen\textsuperscript{a,}\thanks{Corresponding author: robinchen@swu.edu.cn}}
\author[b]{Ren-qian Zhang}

\affil[a]{College of Economics and Management, Research Institute of Intelligent Finance and Platform Economics, Southwest University, Chongqing 400715, China}
\affil[b]{School of Economics and Management, Beihang University, Beijing 100191, China}
\date{}

\begin{document}

\maketitle

\begin{abstract}

This paper investigates a multi-product stochastic inventory problem in which a cash-constrained online retailer can adopt order-based loan provided by some Chinese e-commerce platforms to speed up its cash recovery for deferred revenue.  We first build deterministic models for the problem and then develop the corresponding stochastic programming models to maximize the retailers' expected profit over the planning horizon. The uncertainty of customer demand is represented by scenario trees, and a scenario reduction technique is used to solve the problem when the scenario trees are too large. We conduct numerical tests based on real data
crawling from an online store.  The results show that the stochastic
model outperforms the deterministic model, especially when the retailer is less cash-constrained. Moreover, the retailer tends to choose using order-based loan when its initial available cash is small or facing long receipt delay length.
\end{abstract}


\section{Introduction}

Cash availability is very important to small-to-medium-sized enterprises (SMEs), including start-up companies and small retailers. They usually lack sufficient capital to absorb large losses, and it is difficult for them to obtain external financing such as loans or venture capital compared with large companies due to production capacity and sales scale. Research from the Association of Chartered Certified Accountants (ACCA) and the Institute of Management Accountants (IMA) \citep{Survey2014} found that access to growth capital was one of the two main factors affecting the businesses of SMEs in the UK over the past two years. A report by research firm CB Insights found that 29\% of start-ups failed because of cash crises \citep{cbinsights2018}. Moreover, SMEs tend to be the first to feel the effects of financial crises such as the COVID-19 pandemic \citep{harvard2020}. On Chinese e-commerce platforms such as Taobao.com of Alibaba Group and Jd.com of JD Group, there are many small online retailers operated by only a few people or even a single individual. Generally, the maximum number of products they can purchase depends on their available cash in the current period, which
are similar to nanostores in emerging markets \citep{boulaksil2018cash}.

Many forms of supply chain financing services appear in the market to help SMEs alleviate the cash shortage problem. Thomas Olsen of Bain, a consulting firm, reckons that the supply chain finance market is expanding by 15--25\% per year in the Americas and by 30--50\% in Asia \citep{Eco2018}. Some Chinese e-commerce platforms have also started supply chain financing services for their online retailers in recent years. Alibaba Group opened its own online bank, mybank.cn  in 2010, which mainly provides loans for online stores on its three e-commerce websites: alibaba.com, tmall.com and taobao.com. Its main competitor in China---jd.com---also set up a supply chain finance department in 2012 and provides similar financing services.

\begin{figure}[!ht]
\centering
    \includegraphics[scale=0.85]{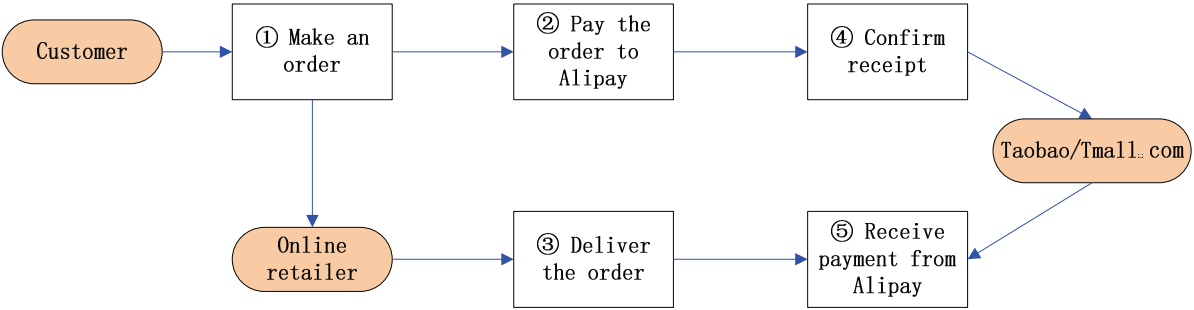}
    \caption{Transaction Process in Taobao.com/Tmall.com.}
    \label{fig:Alipay-process}
\end{figure}

Some Chinese e-commerce platforms provide a financing service called order-based loans that has seldom been addressed in the literature. The aim of this service is to speed up the retailer's deferred revenue, which is common in Chinese online shopping. In the real transactions of Chinese e-commerce platforms, payments from customers are first obtained by the platform. The platform then transfers the payment to the retailer's account after customers confirm the receipt of the order with additional processing time. Therefore, the retailer actually receives the customer's payment after a time lag (usually at least 2 weeks). The transaction process on Taobao.com or Tmall.com is shown in Figure \ref{fig:Alipay-process}, where Alipay is the online payment system of Alibaba Group. The circled numbers in the figure represent the sequence of the transaction process.

Order-based loan is a kind of financing services intended to help speed up retailers' cash recovery. Once it is adopted for an order, the platform immediately makes payments to the retailer on behalf of customers after the order is marked as delivered. The retailer repays the loan with some interest to the platform after an agreed time limit. The retailer does not need to wait for its payment, which is a good method to alleviate the retailer's cash shortage. The payment delay length is zero when applying order-based loan to an order.

Given the above background, the aim of this paper is to investigate the impacts of order-based loan on an online retailer's operational decisions. We make the following contribution to the literature.
\begin{itemize}
\item We consider a supply chain financing service: order-based loan in a stochastic inventory problem, which has been used in business practice but seldom addressed in the inventory literature.

\item A stochastic multi-product inventory model is built for the problem. Scenario trees are adopted to represent the demand uncertainty, and we use a scenario tree reduction method to solve the problem if the tree size is large.

\item By crawling real transaction data from an online store and performing numerical analysis, we find that factors such as the initial available cash, receipt delay length, overhead cost and demand fluctuations may influence a retailer's decision to apply for order-based loans.
\end{itemize}

\section{Literature Review}

The literature associated with our work can be classified into two main streams: inventory management with financing decisions and scenario programming in operational research.

\subsection{Inventory management with financing decisions}

Financial flow management is an important part of supply chain management \citep{cooper1997supply}. In recent years, there have already been many works addressing joint inventory management and financing decisions. Here, we only review some of them for brevity.

\cite{buzacott2004inventory} first considered operational and financial decisions by analyzing a cash-constrained retailer with asset-based financing provided by a bank. \cite{dada2008financing} built a Stackelberg gaming model involving a bank, manufacturer and cash-constrained retailer. \cite{raghavan2011short} considered a problem with a retailer and manufacturer facing cash constraints. \cite{kouvelis2012financing} investigated the structure of optimal trade credit contracts between a supplier and bank. \cite{moussawi2013joint} developed a model considering delayed payments for a cash-constrained retailer. \cite{tunca2018buyer} discussed the role and efficiency of buyer intermediation in supplier financing through a game-theoretic approach. \cite{wu2019trade} proposed a trade credit model with asymmetric competing retailers in which one weak retailer was capital-constrained. \cite{jin2019non} discussed noncollaborative and collaborative financing. \cite{xu2020partial} investigated a supply chain financing system with one supplier and one emission-dependent and capital-constrained manufacturer. \cite{Yuan2020} integrated supply risk into a cash constraint problem. \cite{zhou2020guarantor} investigated manufacturer guarantor financing and third‐party logistics guarantor financing in a four‐party supply chain game.

The abovementioned works mainly built single-period inventory models to analyze problems. In terms of the multi-period setting, \cite{chao2008dynamic} investigated a multi-period self-financing newsvendor problem and proved that the optimal ordering pattern is an approximate base stock policy. \cite{gong2014dynamic} further considered short-term financing. \cite{Katehakis2016Cash} showed that the optimal ordering policy is characterized by a pair of threshold parameters for this problem. \cite{boulaksil2018cash} formulated a cash-constrained stochastic inventory model with consumer loans and supplier credits and obtained managerial insights by simulating numerical cases. \cite{chen2018capital} incorporated cash flow constraints into a multi-period lot-sizing problem with trade credit. \cite{ChenZhang2019} extended this problem with credit-based loans. \cite{bi2018dynamic} addressed short-term financing in an inventory problem with multi-item products. \cite{li2020two} considered a two-stage inventory model with financing and demand updating using a Bayesian approach. \cite{chen2020dynamic} developed an $(s, C(x), S)$ policy for a stochastic inventory problem with fixed costs and cash constraints.

\subsection{Scenario programming in operational research}

In decision making under uncertainty, one method is scenario programming in which the uncertainty is represented by scenario trees. \cite{hoyland2001generating} presented a method to generate a limited number of scenarios that satisfy specified statistical properties. \cite{haugen2001progressive} adopted a progressive hedging algorithm for each scenario to solve a stochastic lot-sizing problem. \cite{heitsch2003scenario} developed two scenario reduction algorithms: new versions of forward and backward type algorithms. \cite{brandimarte2006multi} proposed a fix-and-relax strategy to solve a scenario programming problem with multi-item capacitated lot sizing. \cite{beraldi2006scenario} also made use of this method to solve a lot-sizing problem with uncertain processing times. \cite{helber2013dynamic} took a scenario approach to solve a capacitated lot-sizing problem with random demand and dynamic safety stocks. \cite{feng2013scenario} used the scenario generation and reduction method for a stochastic power generation expansion planning problem. \cite{hu2016two} and \cite{hu2018multi} adopted similar methods to solve two-stage and multi-stage stochastic problems, respectively. \cite{hnaien2017robust} considered a minmax robust lot-sizing problem with uncertain lead time. \cite{fattahi2018multi} used the forward scenario construction technique and Benders decomposition algorithm to solve a multi-period supply chain network redesign problem.

To the best of our knowledge, there is no work in the literature that considers order-based loans in a multi-period stochastic inventory problem, which inspired us to investigate this topic.

\section{Problem description}

For convenience, some notations adopted in this paper are listed in Table \ref{table:notations}. Other relevant notations will be introduced as needed.
\begin{table}[!ht]
\centering
\setstretch{1.5}
\caption{Some notations adopted in the paper.}\label{table:notations}
\begin{tabular}{ll}
\toprule
\multicolumn{2}{l}{Indices:}\\
\specialrule{0em}{1pt}{1pt}
$t$\qquad     &Index of a period, $t=1,\dots,T$.\\
$n$     &Index of an item, $n=1,2,\dots,N$.\\
\specialrule{0em}{3pt}{3pt}

\multicolumn{2}{l}{Operational parameters:}\\
\specialrule{0em}{1pt}{1pt}
$I_{n,0}$\qquad  &Initial inventory for item $n$.\\
$d_{n,t}$ &Demand for item $n$ in period $t$, which is uncertain \\
%, which is stochastic and its expectation is $\overline{d}_{n,t}$.\\
$v_{n}$ &Unit variable ordering cost of the retailer for item $n$.\\
\specialrule{0em}{3pt}{3pt}
\multicolumn{2}{l}{Financial parameters:}\\
\specialrule{0em}{1pt}{1pt}
$C_{0}$\qquad  &Initial cash at the beginning of planning horizon.\\
$H_t$ &Overhead costs in each period (e.g., wages or rents) \\
$L$ & Receipt delay length. \\
$p_{n}$ & Selling price of item $n$. \\
$v_{n}$ & Unit variable ordering cost of item $n$. \\
$\alpha$ &Discount rate for unreceived revenues at the end of planning horizon.\\
$B_U$ & The amount of order-based loan.\\
$r_o$ &Interest rate of order-based loan.\\
\specialrule{0em}{3pt}{3pt}
\multicolumn{2}{l}{Decision variables:}\\
\specialrule{0em}{1pt}{1pt}
$Q_{n,t}$     &Ordering quantity for item $n$ in period $t$. \\
$g_{n,t}$  & The quantity of order-based loan used for product $n$ in period $t$. \\
$R_{n,t}$    &Retailer's revenue from item $n$ in period $t$.\\
$I_{n,t}$ 	 &End-of-period inventory level for item $n$ in period $t$. \\
$C_{t}$ 	 &End-of-period cash balance in period $t$. \\
$FC$ 	 &Final cash at period $T$ after discounting unreceived revenues. \\
\bottomrule
\end{tabular}
\end{table}

Our study is about a cash-constrained online retailer selling multiple products on Tmall.com, which is a Chinese e-commerce platform owned by the Alibaba Group. In the problem, the retailer has initial cash $C_0$ before making operational decisions. There are $N$ products for sale, and its planning horizon length is $T$. The initial inventory for product $n$ at the beginning of the planning horizon is $I_{n,0}$. Overhead costs $H_t$, such as wages or rents, may have to be paid in each period, which are incurred irrespective of the retailer's operation activity.

The retailer purchases products from its suppliers, and the unit ordering cost for product $n$ is $v_n$. The retailer faces uncertain demand from customers. The demand for product $n$ in period $t$ is represented by $d_{n,t}$. In each period, the retailer makes decisions about $Q_{n,t}$, which is the purchasing quantity for product $n$ in period $t$. Excess stock is transferred to the next period as inventory, and unmet demand is lost without paying penalty costs. The end-of-period inventory for product $n$ in period $t$ is $I_{n,t}$. The delivery lead time of products is assumed to be zero. The inventory flow equation for product $n$ for two consecutive periods is given below, where $(x)^+$ denotes $\max\{x,0\}$.

\begin{equation}
  I_{n,t}=\max\{I_{n,t-1}+Q_{n,t}-d_{n,t}, 0\} = (I_{n,t-1}+Q_{n,t}-d_{n,t})^+\label{eq:inventoryFlow-uncertain}
\end{equation}

Assume the time lag length of receipt for product $n$ is $L$. Define $R_{n,t}$ as the realized revenue from product $n$ in period $t$. $R_{n,t}$ can be expressed by Eq. \eqref{eq:revenueUncertain}, where $L$ is the receipt delay length and $\min\{ I_{n,t-L-1}+Q_{n,t-L}, d_{n,t-L}\}$ is the realized demand for item $n$ in period $t$. Based on Eq. \eqref{eq:inventoryFlow-uncertain}, $R_{n,t}$ can also be written as $p_n(I_{n,t-L-1}+Q_{n,t-L}-I_{n,t-L})$ because $\min\{I_{n,t-L-1}+Q_{n,t-L}, d_{n,t-L}\} = I_{n,t-L-1}+Q_{n,t-L} - (I_{n,t-L-1}+Q_{n,t-L}- d_{n,t-L})^+$. Eq. \eqref{eq:revenueUncertain} means that the sales income from period $t$ will be received in period $t+L$.
\begin{equation}
R_{n,t}=\begin{cases}
p_n\min\{I_{n,t-L-1}+Q_{n,t-L}, d_{n,t-L}\}=p_n(I_{n,t-L-1}+Q_{n,t-L}-I_{n,t-L})\quad & t>L\\
0 & t\leq L
\end{cases}\label{eq:revenueUncertain}
\end{equation}

Assume that there are no unreceived payments at the beginning of the planning horizon.  The cash flow equation is:
\begin{equation}
  C_t=\begin{cases}
  C_{t-1}+\sum_{n=1}^{N}R_{n,t}-\sum_{n=1}^{N}v_nQ_{n,t}-H_t& t>0\\
  C_0 &t=0
  \end{cases}\label{eq:cashFlow}
\end{equation}

In Eq. \eqref{eq:cashFlow}, when $t>0$, one period's end-of-period cash balance equals its initial cash in this period ($C_{t-1}$) plus this period's total revenues ($\sum_{n=1}^{N}R_{n,t}$), minus this period's total ordering costs ($\sum_{n=1}^{N}v_nQ_{n,t}$) and overhead costs ($H_t$). At the end of planning horizon $T$, there are still some unreceived payments. For ease of analysis, we assume that the retailer discounts those payments to the final period $T$ with discounting factor $\alpha$. According to the above description, the final cash $FC$ at the end of period $T$ after discounting unreceived payments is shown by Eq. \eqref{eq:finalCashAfterDiscount}, where $C_T$ is the cash balance at period $T$ before discounting the unreceived payments. The total amount of unreceived payments is $\sum_{n=1}^{N}\sum_{k=1}^{L} R_{n,T+k}$.

\begin{equation}
    FC=C_T+\sum_{n=1}^{N}\sum_{k=1}^{L}\frac{1}{(1+\alpha)^{k}} R_{n,T+k}\label{eq:finalCashAfterDiscount}
\end{equation}

The retailer faces cash constraints in each period, which is common in real transactions for such small retailers. The cash constraint means that the retailer's maximum purchasing quantity is bounded by its available cash, which can be represented by the following formula:
\begin{equation}
    \sum_{n=1}^{N}v_nQ_{n,t}\leq C_{t-1}\label{con:cash-constraint}
\end{equation}

\section{Mathematical models}
In this section, several mathematical models are formulated to address different financing solutions provided by e-commerce platforms. We first provide the deterministic models and then build the corresponding stochastic models.

\subsection{Deterministic model}

In the deterministic situation, the uncertain demand for product $n$ in period $t$ is predicted to be the mean value of its demand distribution, which is represented as $\overline{d}_{n,t}$. We build two models below that are distinguished by whether order-based loan is used.

\subsubsection{ Operating with Self-Owned Cash}

For many retailers, if operating without supply chain financing, the cash inflow is mainly initially owned cash and revenues from sales. The whole linear deterministic model for the situation of operating with self-owned cash (SO-D) is built below.

\clearpage
\noindent{\bf Model SO-D}

\begin{align}
\max\quad &FC=C_T+\sum_{n=1}^{N}\sum_{k=1}^{L}\frac{1}{(1+\alpha)^{k}} R_{n,T+k}\tag{\ref{eq:finalCashAfterDiscount}}\\
\text{s.t.}\quad& \forall n, t\nonumber\\
&\eqref{eq:revenueUncertain}, \eqref{eq:cashFlow}, \eqref{con:cash-constraint}\nonumber\\
&I_{n,t}\leq I_{n,t-1}+Q_{n,t}-\overline{d}_{n,t}+(1-\delta_{n,t})M\label{con:inventory-flow-determin1}\\
&I_{n,t}\geq I_{n,t-1}+Q_{n,t}-\overline{d}_{n,t}-(1-\delta_{n,t})M\label{con:inventory-flow-determin2}\\
&I_{n,t-1}+Q_{n,t}-\overline{d}_{n,t}\geq -(1-\delta_{n,t})M\\
&I_{n,t-1}+Q_{n,t}-\overline{d}_{n,t}\leq \delta_{n,t}M\\
%&I_{n,t-1}+Q_{n,t}\geq \overline{d}_{n,t} - (1-\delta_{n,t})M\\
&I_{n,t}\leq \delta_{n,t}M\label{con:inventory-flow-determin3}\\
&I_{n,t}\geq 0, Q_{n,t}\geq 0, \delta_{n,t}\in\{0, 1\}\label{con:BinaryNonNegative-determin1}
\end{align}

The objective function \eqref{eq:finalCashAfterDiscount} is to maximize the expected final cash after discounting the unreceived payments. Constraint \eqref{eq:revenueUncertain} is the revenue equation, while Constraint \eqref{eq:cashFlow} is the cash flow equation. Constraint \eqref{con:cash-constraint} is the cash constraint. Constraints \eqref{con:inventory-flow-determin1}--\eqref{con:inventory-flow-determin3} are the linear expressions for the inventory flow equation \eqref{eq:inventoryFlow-uncertain}. Constraint \eqref{con:BinaryNonNegative-determin1} shows the nonnegativity/binarity of
end-of-period inventory, ordering quantity and $\delta_t$. The nonnegativity of $I_{t}$ is the result of the lost sales assumption.

\subsubsection{Operating with order-based loan}

In the presence of order-based loan, when the retailer's application for order-based loan is approved, the platform sets a credit limit for the retailer. Assume that the limit of order-based loans is $B_U$. The retailer can use this kind of loan repeatedly for many orders over a time limit as long as the total payment subject to order-based loans is less than $B_U$.

The interest rate of the order-based loan is $R_o$. A new decision variable $g_{n,t}$ represents the quantity of order-based loans employed for product $n$ in period $t$. For ease of analysis, assume that the retailer repays an order-based loan after a receipt delay of length $L$, and the time limit for the order-based loan is the planning horizon length $T$. The deterministic model for order-based loans (OL-D) is given below:

\clearpage
\vspace{5pt}
\noindent{\bf Model OL-D}
\vspace{5pt}
\begin{align}
\max\quad &FC=C_T+\sum_{n=1}^{N}\sum_{k=1}^{L}\frac{1}{(1+\alpha)^{k}} R_{n,T+k}\tag{\ref{eq:finalCashAfterDiscount}}\\
\text{s.t.}\quad& \forall n, t\nonumber\\
&\eqref{eq:cashFlow},\eqref{con:cash-constraint}-\eqref{con:BinaryNonNegative-determin1} \nonumber\\
&g_{n,t}\leq I_{n,t-1}+Q_{n,t}-I_{n,t}\label{con:lost-sale-determin}\\
&R_{n,t} =\begin{cases}
p_n(I_{n,t-L-1}+Q_{n,t-L}-I_{n,t-L}-g_{n,t-L}-g_{n,t-L}(1+R_o)^L) & T<t\leq T+L\\
 p_n(g_{n,t}+I_{n,t-L-1}+Q_{n,t-L}-I_{n,t-L}-g_{n,t-L}-g_{n,t-L}(1+R_o)^L) & L<t\leq T\\
  p_n g_{n,t}& t\leq L
\end{cases}\label{eq:revenue:credit-determin}\\
&\sum_{n=1}^N\sum_{t=1}^T p_n g_{n,t}\leq B_U\label{con:order-Total-determin}\\
&g_{n,t}\geq 0\label{con:lostsale-above0-determin}
\end{align}

The differences between Model OL-D and Model SO-D are the constraints related to $g_{n,t}$. Constraint \eqref{con:lost-sale-determin} means that the order-based loan employed for a product should be less than its realized demand ($I_{n,t-1}+Q_{n,t}-I_{n,t}$) in this period. Constraint \eqref{eq:revenue:credit-determin} gives the expression for revenue in the presence of order-based loans: when $t\leq L$, the retailer's revenue from product $n$ is $p_ng_{n,t}$ since it can use order-based loan $g_{n, t}$ for this product; when $t>L$, $I_{n,t-L-1}+Q_{n,t-L}-I_{n,t-L}-g_{n,t-L}$ is the quantity of demand that does not use order-based loan in period $t-L$, and this portion of revenue comes at period $t$; and $p_ng_{n,t-L}(1+R_o)^L)$ is the total principle and interest that the retailer needs to repay in period $t$ for period $t-L$'s order-based loans; when $t>T$, the retailer cannot use order-based loan. Note that because of receipt delaying, the retailer receives some revenue when $t>T$. Constraint \eqref{con:order-Total-determin} suggests that the total employed order-based loan over the planning horizon should be no more than the order-based loan limit. Constraint \eqref{con:lostsale-above0-determin} ensures the nonnegativity of $g_{n,t}$.

\subsection{Stochastic model}

In real transactions, customer demand is uncertain, which is represented by demand scenarios in this paper. The demand for product $n$ in period $t$ under scenario $s$ is represented by $d_{n,t}^{s}$, and the probability of scenario $s$ occurring is $\Pr _s$. The set of all scenarios is $S$. Each scenario consists of the discrete values of all demands over the whole planning horizon. In the scenario programming model, the decision variables are $Q_{n,t}^s$, $R^s_{n,t}$, $I^s_{n,t}$, $g^s_{n,t}$, $C^s_{n,t}$, $FC^s$, and $\delta^s_{n,t}$. The stochastic model when the retailer is operating with self-owned cash (SO-S) is formulated below:

\clearpage
\vspace{5pt}
\noindent{\bf Model SO-S}
\vspace{5pt}
\begin{align}
\max\quad &\sum_{s\in S}\Pr\nolimits_s FC^s=\sum_{s\in S} \Pr\nolimits_s\left[\sum_{n=1}^{N}\sum_{k=1}^{L}\frac{1}{(1+\alpha)^{k}} R^s_{n,T+k}\right]\label{obj-sto1}\\
\text{s.t.}\quad& \forall n, t, s\nonumber\\
&I^s_{n,t}\leq I^s_{n,t-1}+Q^s_{n,t}-d^s_{n,t}+(1-\delta^s_{n,t})M\label{con:inventory-flow-sto1}\\
&I^s_{n,t}\geq I^s_{n,t-1}+Q^s_{n,t}-d^s_{n,t}-(1-\delta^s_{n,t})M\label{con:inventory-flow-sto2}\\
&I^s_{n,t-1}+Q^s_{n,t}-d^s_{n,t}\geq -(1-\delta^s_{n,t})M\label{con:inventory-flow-sto3}\\
&I^s_{n,t-1}+Q^s_{n,t}-d^s_{n,t}\leq \delta^s_{n,t}M\label{con:inventory-flow-sto4}\\
&I^s_{n,t}\leq \delta^s_{n,t}M\label{con:inventory-flow-sto5}\\
&R^s_{n,t}=\begin{cases}
p_n(I^s_{n,t-L-1}+Q^s_{n,t-L}-I^s_{n,t-L})\quad & t>L\\
0 & t\leq L
\end{cases}\label{eq:revenue-sto}\\
  &C^s_t=\begin{cases}
  C^s_{t-1}+\sum_{n=1}^{N}R^s_{n,t}-\sum_{n=1}^{N}v_nQ^s_{n,t}-H_t& t>0\\
  C_0 &t=0
  \end{cases}\label{eq:cashFlow-sto}\\
&\sum_{n=1}^{N}v_nQ^s_{n,t}\leq C^s_{t-1}\label{con:cash-constraint-sto1}\\
&\sum_{s'\in J(s,t)}\Pr\nolimits_{s'}Q^s_{n,t}=Q^s_{n,t}\sum_{s'\in J(s,t)}\Pr\nolimits_{s'}\label{con:nonanti-Q}\\
&\sum_{s'\in J(s,t)}\Pr\nolimits_{s'}I^s_{n,t}=I^s_{n,t}\sum_{s'\in J(s,t)}\Pr\nolimits_{s'}\label{con:nonanti-I}\\
&\sum_{s'\in J(s,t)}\Pr\nolimits_{s'}\delta^s_{n,t}=\delta^s_t\sum_{s'\in J(s,t)}\Pr\nolimits_{s'}\label{con:nonanti-delta}\\
&I^s_{n,t}\geq 0, Q^s_{n,t}\geq 0, \delta^s_{n,t}\in\{0, 1\}\label{con:BinaryNonNegative-sto1}
\end{align}

The objective function \eqref{obj-sto1} is to maximize the expected final cash for all scenarios after discounting the unreceived payments. Constraints \eqref{con:inventory-flow-sto1}--\eqref{con:cash-constraint-sto1} and \eqref{con:BinaryNonNegative-sto1} are the scenario-specific expressions for the constraints in Model SO-D. Constraints \eqref{con:nonanti-Q}--\eqref{con:nonanti-delta} are the non-anticipativity constraints in the scenario programming model, where $J_{s,t}$ represents the sets of scenarios that share the same history as scenario $s$ before period $t$. Non-anticipativity means that all the scenarios with the same history before a given period should result in the same values of decision variables until that period. There are no non-anticipativity constraints for $R^s_{n,t}$ and $C^s_{n,t}$ because they are basically the expressions of $I^s_{n,t}$, $Q^s_{n,t}$ and $\delta^s_{n,t}$.

In a similar way, the stochastic model when the retailer is operating with order-based loan (OL-S) is formulated below:

\noindent{\bf Model OL-S}
\begin{align}
\max\quad &\sum_{s\in S}\Pr\nolimits_s FC^s=\sum_{s\in S} \Pr\nolimits_s\left[\sum_{n=1}^{N}\sum_{k=1}^{L}\frac{1}{(1+\alpha)^{k}} R^s_{n,T+k}\right]\tag{\ref{obj-sto1}}\\
\text{s.t.}\quad& \forall n, t, s\nonumber\\
&\eqref{con:inventory-flow-sto1}-\eqref{con:inventory-flow-sto5}, \eqref{eq:cashFlow-sto}-\eqref{con:BinaryNonNegative-sto1} \nonumber\\
&g^s_{n,t}\leq d^s_{n,t}-w^s_{n,t}\label{con:lost-sale-sto}\\
&R^s_{n,t} =\begin{cases}
p_n(I^s_{n,t-L-1}+Q^s_{n,t-L}-I^s_{n,t-L}-g^s_{n,t-L}-g^s_{n,t-L}(1+R_o)^L) & T<t\leq T+L\\
 p_n(g^s_{n,t}+I^s_{n,t-L-1}+Q^s_{n,t-L}-I^s_{n,t-L}-g^s_{n,t-L}-g^s_{n,t-L}(1+R_o)^L) & L<t\leq T\\
  p_n g^s_{n,t}& t\leq L
\end{cases}\label{eq:revenue:credit-sto}\\
&\sum_{n=1}^N\sum_{t=1}^T p_ng^s_{n,t}\leq B_U\label{con:order-Total-sto}\\
&\sum_{s'\in J(s,t)}\Pr\nolimits_{s'}g^s_{n,t}=g^s_{n,t}\sum_{s'\in J(s,t)}\Pr\nolimits_{s'}\label{con:nonanti-g}\\
&g^s_{n,t}\geq 0\label{con:lostsale-above0-sto}
\end{align}

\section{Data source and scenario settings}

In this section, we describe the data used in the numerical cases and the techniques employed for scenario generation and reduction in the study.

\subsection{Data source}

A case study is derived from an online store that sells computer peripherals such as keyboards, mice, headsets, etc. on the Chinese e-commerce platforms Tmall.com and Jd.com. Since it is difficult to obtain internal sales data from retailers, we crawl customer comments from the retailer' online stores. Assume that 10\% of customers leave comments after buying the products and that there is an approximately one week delay between the date they order and the date they comment, so we can estimate the associated demand. Note that although this seems to be a rought estimate for customers' demand, the scenario generation and reduction methods adopted in this paper can be applied to more exact demand values \citep[e.g.][]{feng2013scenario, hu2018multi}.

For ease of analysis, two weeks is deemed a period in the planning horizon. Because demand typically booms in some months, such as the Single's Day shopping event in November and the semester-opening promotion events in March and September, we classify the planning horizon into two situations: booming-demand periods and normal-demand periods. We fit distribution parameters using a maximum likelihood approach. The empirical study revealed that a log-normal distribution generally provides a good fit for our demand samples. The resulting fitted values of the three main products (keyboard, mouse, headset) are shown in Table \ref{tab: fit-log-normal}; the respective $p$-values from Kolmogorov-Smirnov tests are also reported in the table. Since all the $p$-values are above 0.05, we hold that the demand fits a log-normal distribution at the 95\% confidence level. Because the length of order-based loans provided by the e-commerce platform usually is at most 3 months (6 two-week periods), the length of the planning horizon in our problem is 6 periods.

\begin{table}[!ht]
\centering
\setstretch{1.5}
\caption{Fitted log-normal distribution for the three products.}\label{tab: fit-log-normal}
\begin{tabular}{lllllll}
\toprule
\multirow{3}{*}{Product }  &\multicolumn{6}{c}{Log-normal distribution} \\
&\multicolumn{3}{c}{normal-demand periods} &\multicolumn{3}{c}{ booming-demand periods}\\\cmidrule(lr){2-4}\cmidrule(lr){5-7}
&$\mu$ &$\sigma$  & p-value &$\mu$ &$\sigma$ & p-value\\
\midrule
keyboard     &3.66 &0.60 &0.38 &5.79&0.26&0.58\\
mouse    &4.13 &0.66 &0.62 & 5.91 &0.33 & 0.66\\
headset   &3.54&0.46&0.18 &4.96 &0.18 & 0.96\\
\bottomrule
\end{tabular}
\end{table}

\subsection{Scenario generation}

A limited number of scenarios can be generated to satisfy the specified statistical properties of stochastic distributions. This paper adopts the moment matching method to generate scenarios. The main idea behind this method is to minimize the distances between the specified statistical properties of the generated scenarios \citep{hoyland2001generating}.

We briefly summarize the moment matching technique proposed by \cite{hoyland2001generating}. Recall that the probability of scenario $s$ is $\Pr_s$ and that the set of all scenarios is $S$. Let $\psi$ represent the specified stochastic property of the demand distribution and $\Psi$ be the set of all specified properties. Let $V_\psi$ be the value for the stochastic property $\psi$ and $f(s, \Pr_s)$ be the mathematical expression for $V_\psi$. In the moment matching method, the decision variables are scenario realizations $s$ and the corresponding possibility $\Pr_s$. For example, to generate scenarios for stochastic demand with a normal distribution, the stochastic property $\psi$ can be the mean or variance. $V_\psi$ is the given value for mean or variance. Expression $f(s, \Pr_s)$ is $\sum_{s\in S} s\Pr_s$ for the mean and $\sum_{s\in S} \Pr_s(s-\sum_{s\in S} s\Pr_s)^2$ for the variance. Let $h_\psi$ be the weight for property $\psi$. The mathematical model for the moment matching method is as follows:

\begin{align}
    \min_{s, \Pr_s}\quad  &\sum_{\psi\in\Psi} h_\psi(f(s, \Pr\nolimits_s)-V_\psi)^2\label{obj:scenarioGeneration}\\
    &\sum_{s\in S} \Pr\nolimits_s=1\label{con:psum1}\\
    &\Pr\nolimits_s\geq 0, \quad\forall s\label{con:plarger0}
\end{align}

The objective function \eqref{obj:scenarioGeneration} minimizes the overall weighted squared distance between the given values of specified statistical properties and the corresponding values of mathematical expressions. Constraint \eqref{con:psum1} states that the sum of all demand realization probabilities should be equal to 1. Constraint \eqref{con:plarger0} concerns the nonnegativity of the probability.

To control the size of the scenario tree and avoid possible pitfalls of the moment matching method such as \textit{overspecification} and \textit{underspecification}, \cite{hoyland2001generating} gave a simple formula to find the smallest number of scenarios needed:
\begin{equation}
(D+1)y-1\sim \text{the number of specifications}
\end{equation}

The symbol ``$\sim$" denotes ``close to"; $D$ indicates the dimension of each scenario node vector, and $y$ represents the number of branches from each node. For example, with regard to our three-product, six-period problem, $D$ is 18 for a problem with 3 products and 6 periods. By specifying three moments, the mean, variance and skewness, the number of specifications is $18\times 3=54$. Therefore, by the above formula, we choose to create 3 realizations in each period because $(18+1)\times 3-1\geq 54$.

The fmincon function in MATLAB is adopted to solve this nonlinear optimization model \eqref{obj:scenarioGeneration}--\eqref{con:plarger0}. Multiple starting points are tested, and we select the solution with the smallest objective value sufficiently close to zero. To examine how uncertainty can affect decision making, 3 scenario trees are constructed. The size of each scenario tree is $3^6$. A summary of the demand scenarios for the first period in each scenario tree is shown in Table \ref{tab:scenario-tree}. Recall that since demand is independent in each period, demand realizations and corresponding possibilities are the same for each period. For example, in the first period of scenario tree 1, the probability of demand realizations of 133, 246, and 87 for the three products is 0.1. In the latter periods, the probability of the same realizations is still 0.1.

\begin{table}[!ht]
    \centering
    \caption{Three scenario trees for demands.}\label{tab:scenario-tree}
    \begin{tabular}{llllllll}
    \toprule
    \multicolumn{4}{c}{Normal-demand periods} &\multicolumn{4}{c}{Booming-demand periods}\\
    \cmidrule(lr){1-4}\cmidrule(lr){5-8}
     Probability & Product 1 &Product 2&Product 3&Probability & Product 1 &Product 2&Product 3 \\
      \midrule
         {\it Tree 1} \\
        0.1  & 133  & 246  & 87 &0.286&291 &597 &123\\
        0.598  & 30  & 58  & 39 &0.318 &468&322 &124\\
        0.302 & 49  & 57 & 20  &0.396 &268 &293 &177\\
\specialrule{0em}{3pt}{3pt}
{\it Tree 2} \\
0.102  & 134  & 246  & 88  &0.185 &345 &341 &156\\
        0.694  & 32 & 59  & 37 &0.556 &269&302 &123\\
        0.204  & 54  & 56  & 17  &0.259&481 &611 &184 \\
\specialrule{0em}{3pt}{3pt}
{\it Tree 3} \\
0.103  & 134  & 246 & 87  &0.266 &481 &608 &134\\
        0.476  & 28  & 58  &24&0.34 &317 &311 &181 \\
        0.421 & 46 & 58  & 43 &0.394 &259 &309 &121\\
\bottomrule
    \end{tabular}
\end{table}

\subsection{Scenario reduction}

In each scenario tree, there are $3^6$ scenarios. With a larger planning horizon or more products, the problem may become computationally intractable. For example, if using weekly demand, the length of the planning horizon will be 12, and there will be $3^{12}$ scenarios. Therefore, we also make use of the scenario reduction method proposed by \cite{heitsch2003scenario} to reduce the computational complexity. With regard to the two reduction methods: fast forward selection method (FFS) and backward selection method (BS), we make use of FFS because it outperforms BS when the number of selected scenarios is no more than 25\% of the number of original scenarios\citep{heitsch2003scenario, hu2018multi}.

The main idea behind FFS is described as follows: in the first step, select the scenario with the smallest total weighted distance; then, update the distances and select the scenario among the remaining unselected scenarios; and after selecting enough scenarios as required, add the probabilities of unselected scenarios to the probability of their closest selected scenarios. The pseudocode of FFS is given by Algorithm 1 in the appendix.

\section{Numerical analysis}
In this section, we first discuss the stability test results of the scenario generation and reduction method, then make comparisons between deterministic models and stochastic models, and finally investigate the impacts of different initial cash, overhead cost, receipt delay length, and demand fluctuations on the retailer's choice of using order-based loans.

Selling prices and unit purchasing costs for the three products are listed in Table \ref{tab:price-cost}. Assume that the initial cash balance of the retailer is 20,000. The discount rate $\alpha$ is 1\%, and the interest rate $r_0$ for order-based loan is $1.5\%$. The selling prices and unit purchasing costs of the three products are listed in Table \ref{tab:price-cost}. Other parameter values are as follows: receipt delay length $L$ is 2; order-loan limit $B_u$ is 10,000; overhead cost $H_t$ is 2000 in each period; and the demand patterns are [0, 0, 0, 0, 1, 1], where 0 means normal demand in this period and 1 means booming demand.

\begin{table}[!ht]
\centering
\setstretch{1.5}
\caption{Selling prices and purchasing costs of the three products.}
\begin{tabular}{lllllll}
\toprule
    & keyboard &  mouse & headset\\
    \midrule
Selling price &189 &144 &239\\
Unit purchasing cost&120 & 70& 150\\
\bottomrule
\end{tabular}\label{tab:price-cost}
\end{table}

\subsection{Stability test}

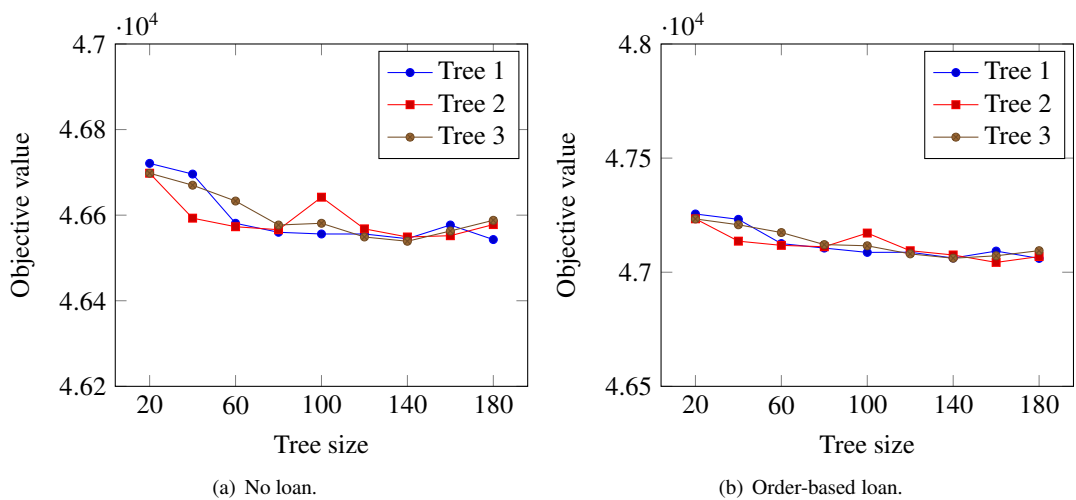
\begin{figure}[!ht]
\centering
\subfigure[No loan.]
{
\begin{tikzpicture}
\begin{axis}
[
xlabel=Tree size, ylabel=Objective value,
scale=0.8, ymin = 46200, ymax = 47000, xtick = {20,60,100, 140,180}, mark size = 1.5pt]
\addplot+[sharp plot] coordinates
{(20,46721)
(40, 46696)
(60,46581)
(80, 46560)
(100,46556)
(120, 46556)
(140,46545)
(160, 46577)
(180,46543)
};
\addplot+[sharp plot] coordinates
{(20,46698)
(40, 46593)
(60,46573)
(80, 46566)
(100,46642)
(120, 46568)
(140,46549)
(160, 46552)
(180,46578)
};
\addplot+[sharp plot] coordinates
{(20,46698)
(40, 46670)
(60,46633)
(80, 46577)
(100,46581)
(120, 46549)
(140,46539)
(160, 46563)
(180,46588)
};
\legend{Tree 1, Tree 2, Tree 3}
\end{axis}
\end{tikzpicture}
}
\subfigure[Order-based loan.]
{
\begin{tikzpicture}
\begin{axis}
[
xlabel=Tree size, ylabel=Objective value,
scale=0.8, ymin = 46500, ymax = 48000, xtick = {20,60,100, 140,180}, mark size = 1.5pt]
\addplot+[sharp plot] coordinates
{
(20, 47255.14646)
(40, 47231.72695)
(60, 47125.36441)
(80, 47105.8101)
(100, 47086.89423)
(120, 47087.54149)
(140, 47062.1865)
(160, 47092.62989)
(180, 47060.55741)
};
\addplot+[sharp plot] coordinates
{
(20,	47234.022)
(40,	47136.11752)
(60,	47117.69456)
(80,	47111.66668)
(100,	47171.55024)
(120,	47094.14907)
(140,	47075.54054)
(160,	47043.03382)
(180,	47068.63297)
};
\addplot+[sharp plot] coordinates
{
(20,	47233.41009)
(40,	47207.95759)
(60,	47173.85551)
(80,	47121.09727)
(100,	47115.87032)
(120,	47080.30799)
(140,	47061.29919)
(160,	47071.98869)
(180,	47094.45884)
};
\legend{Tree 1, Tree 2, Tree 3}
\end{axis}
\end{tikzpicture}
}%
\caption{Scenario sample size stability test.}\label{fig:sample-size-stability}
\end{figure}

We first investigate the relationship between sample size and objective value in Figure \ref{fig:sample-size-stability}. The aims of the test is to balance the scenario size and computational complexity by finding a good scenario sample size such that the objective value is stable.

The scenario sample sizes are generated from 20 to 180 in steps of 20. As Figure \ref{fig:sample-size-stability} shows, there is a slightly decreasing trend in the objective value as the scenario sample size increases. When the scenario size is 120 or 140, the objective value tends to be stable for all the scenario trees in the two situations. Therefore, we decide to select 140 as the scenario sample size. In addition, the objective values have very small fluctuations across different sample sizes.

Next, we analyze the in-sample stability and out-of-sample stability of the scenario tree. In-sample stability means that the objective values should not vary substantially across different scenario trees of the same size, while out-of-sample stability indicates that the solutions from one scenario tree should exhibit close objective values in another scenario tree. Table \ref{tab:stability-test} lists the details of the stability test for the two scenario models. For instance, 46092 in row 1 and column 2 means that substituing the solution from Tree 1 into Tree 2 gets the value 46092.

\begin{table}[!ht]
\centering
\setstretch{1.5}
\caption{Stability test.}
\begin{tabular}{lllllll}
\toprule
&\multicolumn{3}{c}{Model SO-S (no loan)} &\multicolumn{3}{c}{Model OL-S (loan)}\\
    & Tree 1 &  Tree 2 & Tree 3& Tree 1 &  Tree 2 & Tree 3\\
    \midrule
Tree 1   &46067 & 46092 & 46061& 46200 & 46223 & 46193  \\
Tree 2  & 46116 & 46173 & 46024& 46246 & 46306 & 46202  \\
Tree 3  & 46117 & 46126 & 46131& 46247 & 46258 & 46263  \\
\bottomrule
\end{tabular}\label{tab:stability-test}
\end{table}

As Table \ref{tab:stability-test} shows, since the gaps between different solutions in different scenario trees are all less than 5\%, we claim that the stability test for these scenario generation results is valid.

\subsection{Comparison between deterministic model and stochastic model}
By fixing other parameter values and allowing one parameter to vary, we attempt to determine whether it is worth considering demand uncertainty in different parameter settings. The parameter variations are given by Table \ref{table:para-variation}.

\begin{table}[!ht]
\centering
\setstretch{1.5}
\caption{Parameter value variations.}
\begin{tabular}{ll}
\toprule
Parameter & Value variations\\
\midrule
Initial cash $C_0$\qquad~~~~  &10,000\quad 15,000\quad 20.000\quad 25,000\quad 40,000\\
Receipt delay length $L$ &0\quad 1\quad 2\quad 3\quad 4 \\
Overhead Cost $H_t$ &0\quad 1,000\quad 2,000\quad 3,000\quad 4,000\\
\bottomrule
\end{tabular}\label{table:para-variation}
\end{table}

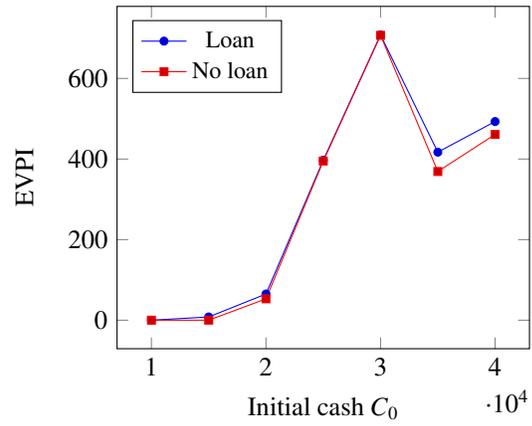
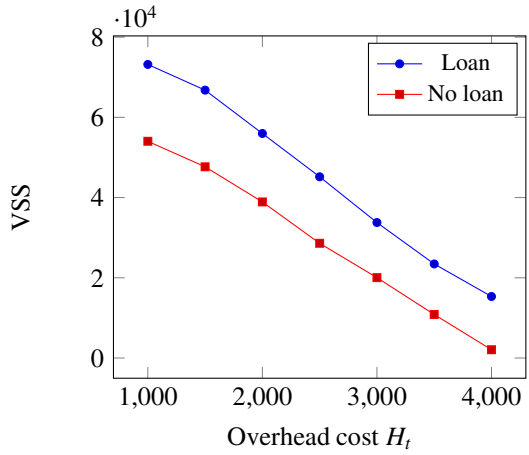
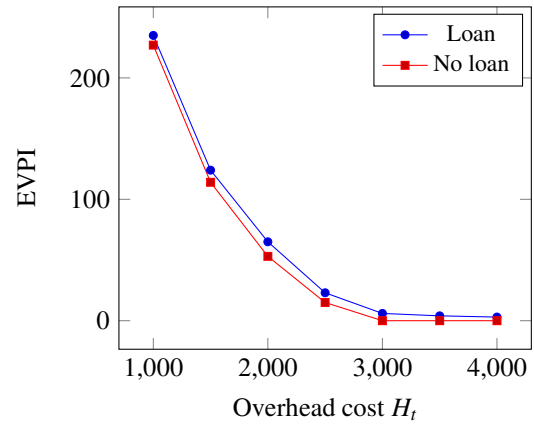
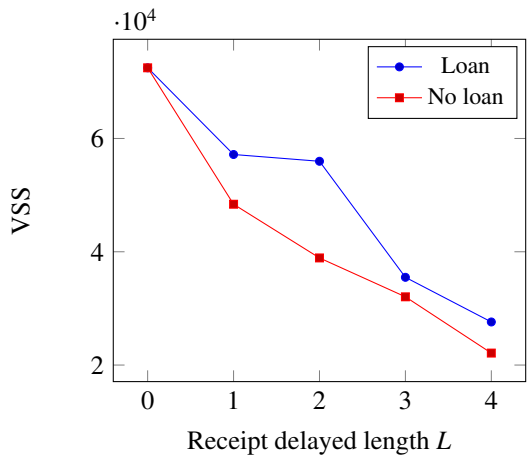
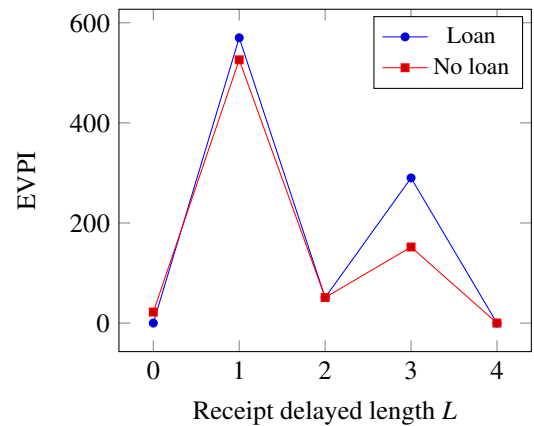
\begin{figure}
\centering
\subfigure[VSS for different initial cash.]
{
\begin{tikzpicture}
\pgfplotsset{every axis legend/.append style={
at={(0.4,0.96)}}, legend style={font=\small}}
\begin{axis}
[%yticklabel=\pgfmathprintnumber{\tick}\,$\%$,
xlabel=Initial cash $C_0$,ylabel=VSS,
scale=0.8, mark size = 1.5pt]
\addplot+[sharp plot] coordinates
{(10000,5673)
(15000,26584)
(20000,55960)
(25000,72824)
(30000,63475)
(35000,63198)
(40000,75824)
};
\addplot+[sharp plot] coordinates
{(10000,29)
(15000,16430)
(20000,36388)
(25000,48529)
(30000,54327)
(35000,58639)
(40000,63358)
};
\legend{Loan, No loan}
\end{axis}
\end{tikzpicture}
}%
\subfigure[EVPI for different initial cash.]
{
\begin{tikzpicture}
\pgfplotsset{every axis legend/.append style={
at={(0.4,0.96)}}, legend style={font=\small}}
\begin{axis}
[
xlabel=Initial cash $C_0$, ylabel=EVPI,
scale=0.8, mark size = 1.5pt]
\addplot+[sharp plot] coordinates
{(10000,0)
(15000,8)
(20000,65)
(25000,397)
(30000,708)
(35000,417)
(40000,493)
};
\addplot+[sharp plot] coordinates
{(10000,0)
(15000,0)
(20000,53)
(25000,395)
(30000,707)
(35000,369)
(40000,461)
};
\legend{Loan, No loan}
\end{axis}
\end{tikzpicture}
}
\subfigure[VSS for different overhead cost.]
{
\pgfplotsset{legend style={font=\small}}
\begin{tikzpicture}
\begin{axis}
[
xlabel=Overhead cost $H_t$,ylabel=VSS, scale=0.8, mark size = 1.5pt ]
\addplot+[sharp plot] coordinates
{(1000,73187)
(1500,66781)
(2000,55974)
(2500,45178)
(3000,33773)
(3500,23450)
(4000,15339)
};
\addplot+[sharp plot] coordinates
{(1000,54027)
(1500,47677)
(2000,38911)
(2500,28601)
(3000,20043)
(3500,10852)
(4000,2073)
};
\legend{Loan, No loan}
\end{axis}
\end{tikzpicture}
}
\subfigure[EVPI for different overhead cost.]
{
\pgfplotsset{legend style={font=\small}}
\begin{tikzpicture}
\begin{axis}
[
xlabel=Overhead cost $H_t$,ylabel=EVPI,scale=0.8, mark size = 1.5pt ]
\addplot+[sharp plot] coordinates
{(1000,235)
(1500,124)
(2000,65)
(2500,23)
(3000,6)
(3500,4)
(4000,3)
};
\addplot+[sharp plot] coordinates
{(1000,227)
(1500,114)
(2000,53)
(2500,15)
(3000,0)
(3500,0)
(4000,0)
};
\legend{Loan, No loan}
\end{axis}
\end{tikzpicture}
}
\subfigure[VSS for different delayed length.]
{
\pgfplotsset{legend style={font=\small}}
\begin{tikzpicture}
\begin{axis}
[
xlabel=Receipt delayed length $L$,ylabel=VSS,scale=0.8, mark size = 1.5pt ]
\addplot+[sharp plot] coordinates
{(0,72475)
(1,57170)
(2,55974)
(3,35489)
(4,27605)
};
\addplot+[sharp plot] coordinates
{(0,72475)
(1,48388)
(2,38911)
(3,32047)
(4,22113)
};
\legend{Loan, No loan}
\end{axis}
\end{tikzpicture}
}
\subfigure[EVPI for different delayed length.]
{
\pgfplotsset{legend style={font=\small}}
\begin{tikzpicture}
\begin{axis}
[%yticklabel=\pgfmathprintnumber{\tick}\,$\%$,
xlabel=Receipt delayed length $L$,ylabel=EVPI,scale=0.8, mark size = 1.5pt ]
\addplot+[sharp plot] coordinates
{(0,0)
(1,570)
(2,51)
(3,290)
(4,0)
};
\addplot+[sharp plot] coordinates
{(0,22)
(1,526)
(2,51)
(3,152)
(4,0)
};
\legend{Loan, No loan}
\end{axis}
\end{tikzpicture}
}
\caption{Comparison between the stochastic model and deterministic model.}\label{fig:sto-determin}
\end{figure}

Let DV represent the expected value for the solution  of the deterministic model in the stochastic model, SV represent the final value of the stochastic model and PV represent the final value if the retailer has perfect knowledge of the stochastic demand. The expected value of perfect information (EVPI) measures how much one wants to pay for perfect information, and the value of the stochastic solution (VSS) implies the difference between the deterministic model and stochastic model. Apparently, EVPI = PV - SV and VSS = SV - DV. A summary of the comparisons between the deterministic model and stochastic model is illustrated in Figure \ref{fig:sto-determin}.

In Figure \ref{fig:sto-determin}, VSS grows with increasing initial cash but decreases with increasing overhead cost and receipt delay length. Since a larger overhead cost or receipt delay length would make the retailer more cash-constrained during operation, this means that the stochastic model has a greater advantage when the retailer is less cash-constrained. With regard to EVPI, it grows as the initial cash increases and decreases as the overhead cost increases. However, as the receipt delay length increases, the EVPI does not show a clear monotonous pattern. Moreover, both the VSS and EVPI values when using order-based loans are larger than those when not using loans. Therefore, it is much better to use the stochastic model when using order-based loans.

\subsection{Analysis for using order-based loan}

\begin{figure}[!ht]
\centering
\subfigure[Expected profit gap for different initial cash.]
{
\begin{tikzpicture}
\begin{axis}
[yticklabel=\pgfmathprintnumber{\tick}\,$\%$,
xlabel=Initial cash $C_0$,ylabel=Profit gap,
scale=0.8]
\addplot+[sharp plot] coordinates
{(10000,6.51)
(15000,0.34)
(20000,0.14)
(25000,0.09)
(30000,0.04)
};
\end{axis}
\end{tikzpicture}
}%
\subfigure[Expected profit gap for different receipt delayed length.]
{
\begin{tikzpicture}
\begin{axis}
[yticklabel=\pgfmathprintnumber{\tick}\,$\%$,
xlabel=Receipt delayed length $L$,
scale=0.8]
\addplot+[sharp plot] coordinates
{(0,0)
(1,0.04)
(2,0.14)
(3,6.18)
(4,13.24)
};
\end{axis}
\end{tikzpicture}
}
\subfigure[Expected profit gap for different overhead cost.]
{
\begin{tikzpicture}
\begin{axis}
[yticklabel=\pgfmathprintnumber{\tick}\,$\%$,
xlabel=Overhead cost $H_t$,ylabel=Profit gap, scale=0.8 ]
\addplot+[sharp plot] coordinates
{(1000,0.09)
(1500,0.11)
(2000,0.14)
(2500,0.18)
(3000,0.26)
(3500,0.49)
(4000,0)
(4500,0)
};
\end{axis}
\end{tikzpicture}
}
\subfigure[Expected profit gap for different demand patterns.]
{\label{fig:profit-gap-demand}
\begin{tikzpicture}
\begin{axis}
[yticklabel=\pgfmathprintnumber{\tick}\,$\%$, xtick = {1, 2, 3, 4},
xlabel=demand pattern,scale=0.8 ]
\addplot+[sharp plot] coordinates
{(1,0)
(2,0.14)
(3,1.03)
(4,1.02)
};
\end{axis}
\end{tikzpicture}
}
\caption{Profit gaps with respect to whether order-based loans are used.}\label{fig:profit-gap}
\end{figure}
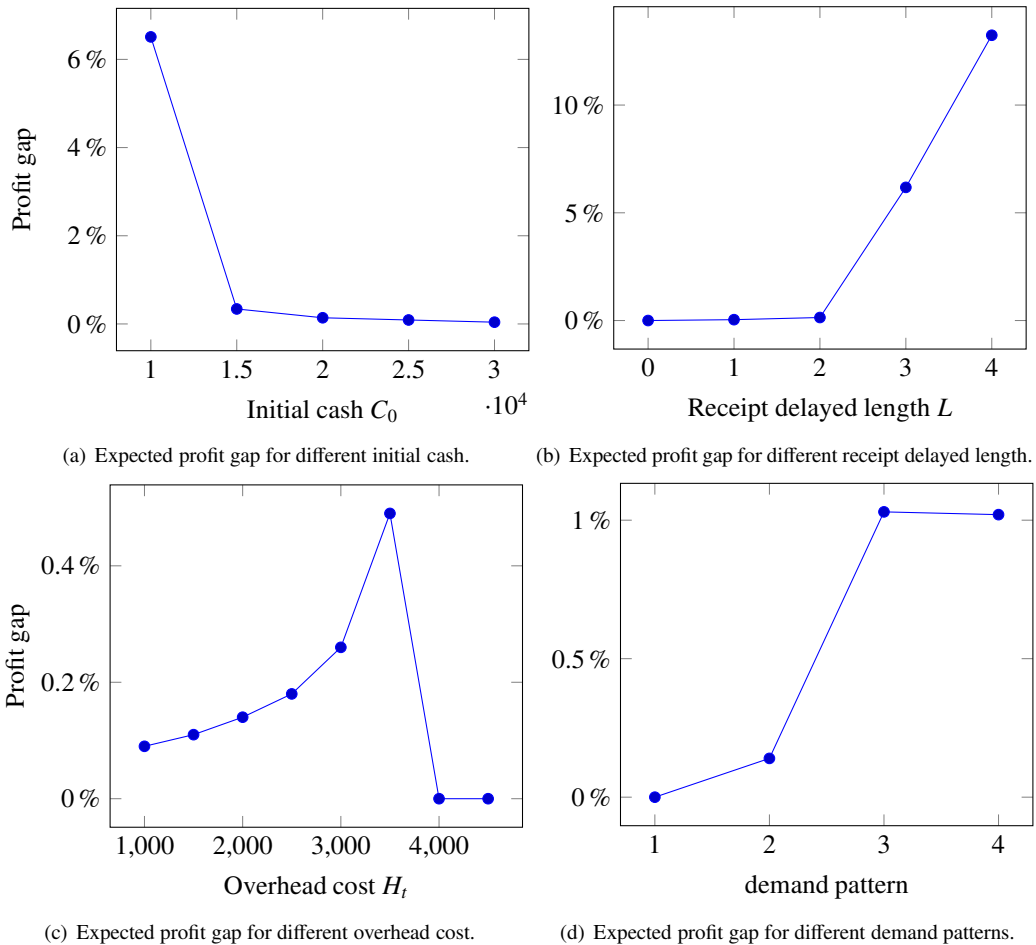

We also want to know which factors affect the retailer's choice of using order-based loans, namely, in which parameter settings the order-based loan model (Model OL-S) achieves higher expected profit than the model of not using loan (SO-S). This is shown by Figure \ref{fig:profit-gap}, where the y axis is the expected profit gap between using order-based loan and not using loan. We also compare several demand patterns in Figure \ref{fig:profit-gap-demand}. The four demand patterns are [0, 0, 0, 0, 0, 0], [0, 0, 0, 0, 1, 1], [1, 1, 0, 0, 0, 0], and [1, 1, 1, 1, 1, 1]. Recall that 0 represents normal demand for this period, while 1 represents booming demand.

The above figures show that the retailer would prefer using order-based loans when the receipt delay length is long or its initial cash is small. When the overhead cost is large, using order-based loans tends to yield more profits than not using them. However, when the overhead cost is too large, there is no advantage in using order-based loans because the retailer cannot make profits. With regard to demand patterns, it seems that it is better for the retailer to use order-based loans when there expected demand is high.

\section{Conclusions}

Many small online retailers encounter cash shortages during operation, and some Chinese e-commerce platforms have provided financing services to help retailers alleviate this problem. In this paper, we discuss order-based loans by building a multi-period multi-product model with demand uncertainty. Scenario programming and reduction techniques are used to solve this problem. Numerical examples show that the retailer tends to use order-based loans when the receipt delay length is long or initial cash is small. When the overhead cost in each period is large but not too large to make the retailer unprofitable, the retailer would also prefer using order-based loans. Future research may consider several directions: one possible direction is to
consider demand substitution in which one product can be substituted by others if it is out of stock; a second possible direction is to investigate other supply chain financing behaviors in a multi-period stochastic problem.

\section*{Acknowledgments}
The research presented in this paper is supported by Chongqing Social Science Planning Fund under Project 2020BS49, the Fundamental Research Funds for the Central Universities of China under Project SWU1909738 and the Research Funds of the Research Institute of Intelligent Finance and Platform Economics, Southwest University under Project 20YJ0105.

\bibliography{liter}
\bibliographystyle{plainnat}

\clearpage
\begin{appendices}
\renewcommand\sectionname{Appendix}

\section{FFS algorithm for scenario reduction}

$s_i$ and $s_j$ represent scenario $i$ and scenario $j$ ($i, j\in S$). $d_{i,j}^{[k]}$ is the distance of $s_i$ and $s_j$ in step $k$ and $WD_{i}^{[k]}$ is the overall weighted distance of scenario $i$ in step $k$. The target number of selected scenarios is $K$. In step $k$, the set of selected scenarios is represented by $\Omega^{[k]}$ while the set of unselected scenarios is represented by $J^{[k]}$. $L(i)$ represents the set of scenarios whose closest scenario is scenario $i$ ($i\in\Omega^{[K]}$) in the final step.

\begin{algorithm}
\setstretch{1.35}
\caption{Fast Forward Algorithm (FFS)}\label{algorithm1}
\KwData{the set of all scenarios $S$, target number of selected scenarios $K$. }
\KwResult{the selected $K$ scenarios. }
Initialize:~~$k\leftarrow 1$,~~$J^{[1]}\leftarrow S$,~~ $\Omega^{[1]}\leftarrow \emptyset$\;

\While {$k<K$}{
\uIf{$k=1$}{
$d^{[1]}_{i,j}$ is the Euclidean distance of scenario $i$ and $j$, $\forall i,j\in J^{[1]}$\;
$WD_i^{[1]}\leftarrow \sum_{s_j\in J^{[1]}}\Pr(s_j)d^{[1]}_{j,i}, \forall i\in J^{[1]}$\;
$ l\leftarrow \arg\min WD_i^{[1]}$\;
$J^{[1]}\leftarrow J^{[1]} \setminus s_{l}$,\quad $\Omega^{[1]}\leftarrow \Omega^{[1]} \cup s_{l}$\;
}
 \Else{
 $d_{j,i}^{[k]}\leftarrow \min\{d_{j,i}^{[k-1]}, d_{j,l}^{[k-1]}\}, \forall i,j\in J^{[k]}$\;
 $WD_i^{[k]}\leftarrow \sum_{s_j\in J^{[k]}}\Pr(s_j)d_{j,i}^{[k]}, \forall i\in J^{[k]}$\;
 $l\leftarrow \arg\min WD_i^{[k]}$\; $J^{[k]}\leftarrow J^{[k]} \setminus s_{l}$,\quad  $\Omega^{[k]}\leftarrow \Omega^{[k]} \cup s_{l}$\;
 }
 $k\leftarrow k+1$\;
}
$\Pr(i)\leftarrow \Pr(i) + \sum_{j\in L(i)}\Pr(j), \forall i\in \Omega^{[K]}$,  $\forall j\in J^{[K]}$, where $L(i)=\{s_j\mid \arg\min_{i\in\Omega^{[K]}, j\in J^{[K]}}\{ d^{[K]}_{i, j}\}=i\}$.
\end{algorithm}

\end{appendices}

\end{document}